\documentclass{amsart}
\usepackage{amssymb}

\usepackage{epic, eepic}

\def\squarebox#1{\hbox to #1{\hfill\vbox to #1{\vfill}}} 
\newcommand{\stopthm}{\hfill\hfill\vbox{\hrule\hbox{\vrule\squarebox 
                 {.667em}\vrule}\hrule}\smallskip}

\newcommand{\Op}{{\operatorname{Op}^{{w}}_h}}
 
\newcommand{\F}{{\mathcal F}} 
\newcommand{\CC}{{\mathbb C}}
\newcommand{\CI}{{\mathcal C}^\infty }
\newcommand{\CIc}{{\mathcal C}^\infty_{\rm{c}} }

\newcommand{\ZZ}{{\mathbb Z}}

\newcommand{\RR}{{\mathbb R}}

\newcommand{\NE}{n_{\rm{e}}}
\newcommand{\NR}{n_{\rm{rh}}}
\newcommand{\NC}{n_{\rm{ch}}}
\newcommand{\EE}{{\rm{e}}}
\newcommand{\RH}{{\rm{rh}}}
\newcommand{\CH}{{\rm{ch}}}

\newcommand{\NN}{{\mathbb N}}

\newcommand{\tr}{\operatorname{tr}}

\renewcommand{\Re}{\mathop{\rm Re}\nolimits}
\renewcommand{\Im}{\mathop{\rm Im}\nolimits}

\theoremstyle{plain}

\newtheorem{thm}{Theorem}

\newtheorem{prop}{Proposition}[section]

\newtheorem{coro}[prop]{Corollary}

\newtheorem{lem}[prop]{Lemma}

\theoremstyle{definition}

\numberwithin{equation}{section}

\title[Normal forms in inverse problems]{Birkhoff normal forms in semi-classical inverse problems}
\author{A. Iantchenko}
\address{Malm{\"o} H{\"o}gskola, Teknik och Samh{\"a}lle, SE-20506, Malm{\"o}, Sweden}
\email{ai@ts.mah.se}
\author{J. Sj{\"o}strand}
\address{ Centre de Math{\'e}matiques, {\'E}cole Polytechnique, FR-91128, Palaiseau
Cedex, France, and UMR 7640 of CNRS}
\email{johannes@math.polytechnique.fr}
\author{M. Zworski}
\address{Department of Mathematics, University of California\\
Berkeley, CA 94720, USA}
\email{zworski@math.berkeley.edu}

\begin{document}
\maketitle

\section{Introduction}
\label{s1}

The purpose of this note is to apply the recent results on 
semi-classical trace formul{\ae} \cite{SjZw},
and on quantum Birkhoff normal forms for semi-classical Fourier
Operators \cite{IS} to inverse problems. We show how 
the {\em classical Birkhoff normal form} can be recovered from
semi-classical spectral invariants. In fact the full 
quantum Birkhoff normal form of the quantum Hamiltonian 
near a closed orbit, and infinitesimally with respect to 
the energy can be recovered. This generalizes recent
results of Guillemin \cite{G} and Zelditch \cite{Zel1},\cite{Zel2},
\cite{Zel} obtained
in the high energy setting (a special case of semi-classical
asymptotics). 

We will illustrate the results in a new setting to which they apply.
Let $ P ( h ) $
be a semi-classical Schr{\"o}dinger operator:
\begin{equation} 
\label{eq:scs}
P ( h )  = - h^2 \Delta + V ( x ) - E  \,, \ \ V \in \CI ( \RR^{n+1} ; \RR) \,,
\end{equation}
with the principal symbol $ p ( x , \xi ) = \xi^2 + V ( x ) - E $. We make
the following assumptions:
\begin{equation}
\label{eq:ass}
\exists \; \epsilon > 0 \,, \  \   p^{-1} ( [-\epsilon , \epsilon] ) \Subset
T^* \RR^{n+1} \,, \ \ dp   \neq 0 \hbox{ on }p^{-1}(0).
\end{equation}
The first assumption guarantees that $ P ( h ) $ has a discrete spectrum near
$ 0 $, and the second one, that the energy surface $ p^{-1} ( 0 ) $ is
smooth.

The semi-classical inverse problem can be formulated as follows:

\begin{quote}
What information about the energy surface, $ \xi^2 + V( x ) = E $, can 
be recovered from the asymptotics of the spectrum of $ - h^2 \Delta 
+ V ( x ) $ in a small fixed neighbourhood of $ E $, as $ h \rightarrow 
0 $\, ?
\end{quote}

This is  very natural and  can be considered as a mathematical
formulation of the general problem arising in 
spectroscopy. In concrete physical situations, the relation
between the Hamiltonian and the spectrum is investigated using 
{\em spectroscopic Hamiltonians} which are essentially Hamiltonians
in Birkhoff normal forms -- see \cite{SR} for one of the first uses
of Birkhoff normal forms in theoretical chemistry, and \cite{JJTF} for
recent applications using experimental spectral data.

The energy surface, $ p^{-1} ( 0 ) \subset T^* \RR^n $ has a natural
bicharactistic foliation, coming from the Hamilton vector field,
$ H_p = \sum_{ j = 1}^n \partial_{\xi_j} p \; \partial_{x_j } - 
\partial_{x_j} p \; \partial_{\xi_j } $. Since it describes classical 
dynamics, its properties are a natural target of a 
recovery procedure based on quantum information, such as the spectrum. 
A na{{\"\i}}ve intuition suggests that the quantum bound states (corresponding
to eigenvalues) should correspond to the closed orbits of the classical
dynamics. The simplest symplectic invariants associated to such closed orbits,
$ \gamma $, are the actions,
\[  S_\gamma = \int_{ \gamma} \xi dx \,, \]
the eigenvalues of the linear Poincar{\'e} map, $ P_\gamma $, and the 
{\em Maslov index}, $ \nu_\gamma \in \ZZ_4 $. We recall that $ P_\gamma $ is
the differential of the (non-linear) Poincar{\'e} map, $ \kappa_\gamma $, 
defined as follows: let $ U $ be a hypersurface in $ p^{-1} ( 0 ) $, 
transversal to $ \gamma $, locally defined near some point on $ \gamma $.
Then  $ U $ is a symplectic manifold. The symplectic map $ \kappa_\gamma :
U \rightarrow U $ is defined as the {\em first return map}:
\[ \kappa_\gamma ( m ) = \exp ( T(m) H_p ) ( m ) \in U \,, \ \
\exp ( t H_p ) ( m ) \notin U \,, \ 0 < t < T( m ) \,, \]
where, strictly speaking, the definition may require replacing $ U $ by 
a larger neighbourhood as the range of $ \kappa_\gamma $.

The relation between the spectrum and closed orbits is provided
by trace formul{\ae}: the Gutzwiller and Balian-Bloch trace formul{\ae}
in physics 
and their mathematical precursor, the Selberg trace formula, and
successors, the Duistermaat-Guillemin and the Guillemin-Melrose formul{\ae}
(see \cite{SjZw} and references given there). Gutzwiller's
formula states that if 
$ \det ( I - P_\gamma^{k} ) \neq 0 $, for $0\neq kT_\gamma \in {\rm
  supp\,}g$,  then 
for $ g , \chi \in \CIc ( \RR ) $, $ \chi \equiv 1 $ near $ 0 $, $ g \equiv 
0 $ near $ 0 $,
\begin{equation}
\label{eq:trace1}
\tr \; \hat g \left( \frac{P ( h ) }{ h } \right) \chi ( P ( h )) 
\sim  \sum_{ \gamma} \sum_{ 0\ne k = - \infty}^\infty
\frac{ e^{ i k S_\gamma/ h + i  \nu_{\gamma, k } \frac{\pi}{2}   }
\;   T_\gamma
 }{ | \det ( P_\gamma^k - I ) |^{\frac12} } 
\sum_{ j=0 }^\infty a_{j , \gamma, k 
 } ( g)  h^j  \,, \ \  a_{ 0, \gamma , k } =   g  ( k  T_\gamma ) \,. 
\end{equation}
where $ T _\gamma $ is the period of $ \gamma $, and the $ \gamma $
summation is over all simple closed orbits of $ H_p $ on $ p^{-1} ( 0 ) $.
The Maslov index, $ \nu_{\gamma, k } $, may in principle depend on the
number of times of ``going around'' the orbit. The coefficients $ 
a_{ j , \gamma, k } $ are distributions on $ \RR $: $ \CIc ( \RR ) 
\ni g \mapsto a_{ j , \gamma , k } $.

It is now standard 
that under the non-degeneracy and simplicity assumptions, 
$ \det ( I - P_\gamma ) \neq 0 $, 
$ T_\gamma =
T_{\gamma'} \; \Rightarrow \; \gamma =  \gamma' $,  
we can recover the actions, $ S_\gamma $,  periods, $ T_\gamma $, 
Maslov indices, $ \nu_{\gamma, k } $, and $ |\det ( I - P_\gamma^k ) |$, 
for all $ k $. 
Under stronger non-degeneracy conditions 
(see \eqref{eq:dio1} and \eqref{eq:dio2} below), 
a theorem of Fried \cite{F} (see Proposition \ref{p:3.2}  below) then 
shows that the eigenvalues of the linear Poincar{\'e} map, $ P_\gamma $, can be 
recovered. 

The higher order coefficients $ a_{ j ,\gamma , k } $ contain further
information and they form a family of semi-classical spectral invariants. 
It turns out that under a stronger version of the non-degeneracy 
hypothesis, we can recover from them the full infinitesimal 
information about the non-linear Poincar{\'e} map $ \kappa_\gamma $. 
For simplicity we state it here in the case of elliptic orbits only,
which is the semi-classical analogue of Guillemin's result \cite{G} -- 
see Sect.\ref{s5} for the general statement. 

Thus suppose that 
the eigenvalues of 
$ d \kappa_\gamma ( w_0 ) $, $ w_0 = U \cap \gamma $, 
are given by $ \exp ( \pm i \theta_j ) $, 
$ j = 1, \cdots, n $, $ \theta_j \in ( 0 ,  \pi ) $. We strengthen the
non-degeneracy assumption to independence from $ 2 \pi $ over rationals:
\begin{equation}
\label{eq:dio1}
\sum_{ j = 1}^{n} k_j \theta_j \in 2 \pi \ZZ 
\,, \ \ k_j \in \ZZ
\ \Longrightarrow 
k_j = 0 \,, \ j = 1 , \cdots , n \,.
\end{equation}

The now classical normal form theorem of Birkhoff, Lewis and Sternberg
says that
there exist symplectic coordinates $ ( x, \xi ) \in \RR^{2n} $,
centered at $ w_0 $, such that in a neighbourhood of $ w_0 $,
we have 
\begin{gather}
\label{eq:ellbnf}
\begin{gathered}
\kappa_\gamma  = \kappa_0 + \kappa^\flat \,, \ \ \kappa^\flat
 ( x , \xi ) = {\mathcal O} ( ( x^2 + \xi^2 )
^\infty ) \,, \\
\kappa_0 ( \iota , \omega ) = \left(  \iota , \omega - \frac{\partial p}{
\partial \iota } \right) \,, \ \ 
\iota = ( \iota_1, \cdots , \iota_n ) \,, \ \omega = ( \omega_1, \cdots , 
\omega_n ) \,, \\
\iota_j = \frac{1}{2}( x_j^2 + \xi_j^2 ) \,, \ \ \omega_j = 
\arg
 ( x_j + i \xi _j ) \,, \\
p ( \iota) = p ( 0 ) - \sum_{j=1}^n 
{ \theta_j } \iota_j + {\mathcal O}( \iota^2 ) 
\,, 
\end{gathered}
\end{gather}
see \cite[Theorem 2.1]{G}. The Taylor series of the twist map $ \kappa_0 $
is called the {\em Birkhoff normal form} of $ \kappa $. It invariantly 
describes the infinitesimal properties of the flow near $ \gamma$.
Since we chose $ \theta_j$'s in a unique way, $ 0 < \theta_j < \pi $,
the Birkhoff normal form is unique.

Specialized to this case our general result Theorem \ref{t:4} can 
be stated as
\begin{thm}
\label{t:0}
Suppose that $ P ( h ) = -h^2 \Delta + V ( x ) - E $ is a 
semi-classical Schr{\"o}dinger operator satisfying \eqref{eq:ass},
and that $ \gamma $ is an elliptic closed orbit of $ H_p $ on $ p^{-1} ( 0 ) $,
with eigenvalues of $ P_\gamma $, $ e^{ \pm i \theta_j } $, satisfying 
\eqref{eq:dio1}. Then the leading term and the coefficients
$ a_{ j , \gamma , k }$
in \eqref{eq:trace1},
determine the Birkhoff normal form of
$ {\kappa}_\gamma $.
\end{thm}

As stated above, our motivation for the recovery of the classical 
Birkhoff normal
forms came from an attempt to understand the results of Guillemin \cite{G}
and Zelditch \cite{Zel1},\cite{Zel2},\cite{Zel}, 
in the context of the recent trace formula of Sj{\"o}strand-Zworski \cite{SjZw},
and of the (quantum) Birkhoff normal form for semi-classical Fourier integral 
operator of Iantchenko \cite{I} and Iantchenko-Sj{\"o}strand \cite{IS}
(see also \cite{Sj}).
Guillemin and Zelditch considered the high energy (or $\CI $ singularities)
case, corresponding to $ h^2 \Delta_g - 1 $ semi-classically. Guillemin's
starting point was his earlier work with Fran{\c c}oise \cite{GF},
and an observation of Zelditch on the non-commutative residues of 
(classical) Fourier integral operator \cite{GZ}. Zelditch's approach
also used the non-commutative residue but was more concrete and 
computational. 

In addition to recovering the classical Birkhoff normal form, \cite{G},
\cite{Zel1}, \cite{Zel2}, \cite{Zel} were concerned with recovering 
the full quantum Birkhoff normal form of the operator near the closed 
trajectory. Without giving a definition of the Birkhoff normal form this 
result follows from 
\begin{thm}
Suppose that the assumptions of Theorem \ref{t:0} are satisfied for 
potentials $ V $ and $ \widetilde V $, with orbits $ \gamma $ and 
$ \tilde \gamma $, respectively. If, in the 
notation of \eqref{eq:trace1}, $ S_\gamma  = \tilde S_{\tilde 
\gamma }$, $ T_\gamma = \widetilde T_{\tilde \gamma } $, and 
 $ |\det ( I - P_\gamma^k ) | = |\det ( I - \widetilde
P_{\tilde \gamma}^k )|
  $, $ \nu_{ \gamma, k} = \tilde \nu_{ \tilde \gamma, k}$,
  for all $ k $, then 
\[ \forall \; j, k \,, \ \ a_{ j , \gamma , k } = \tilde a_{ j , \tilde 
\gamma , k } \ \Longleftrightarrow \ \left\{ \begin{array}{l}
\text{There exists a unitary $h$-Fourier integral operator
 $F$,} \\
\text{ microlocally defined 
in a neighbourhood  of $\gamma $,}\\
\text{ such that
$  ( -h^2 \Delta + \widetilde{V}) F
=F ( - h^2 \Delta + V )   $, }\\
\text{microlocally to infinite order along $ \gamma $.} 
\end{array} \right.
\,.\]
\end{thm}
The general result follows from Proposition \ref{p:n.1} and 
Theorem \ref{t:4}
below, and is given in Corollary \ref{c:n}. The 
definitions of $h$-Fourier integral operators, and of 
agreement ``to infinite order'' are recalled in 
Sect. \ref{s2}.

Our approach is a straightforward application of \eqref{eq:trace2} and
\eqref{eq:qbnf} quoted below from our earlier papers, and of a result 
of Fried \cite{F} on the recovery of the eigenvalues of 
$P_\gamma$ (see Proposition \ref{p:3.2} below). The trace formula 
shows that the traces of powers of 
the quantum monodromy operator of \cite{SjZw} are spectral invariants.
The quantum monodromy operator is a semi-classical operator which 
quantizes the Poincar{\'e} map for a given orbit $ \gamma $, and once 
it is put into the Birkhoff normal form using \eqref{eq:qbnf} the traces
of its powers are easy to compute -- see Proposition \ref{p:3.1}. 
The recovery of the classical Birkhoff normal form is then very clear -- see 
Theorem \ref{t:1}. 

One of the most striking applications of this approach to inverse 
problems is the result
of Zelditch \cite{Zel} on the recovery of bi-axial analytic planar
domains from the spectrum of the (Dirichlet or Neumann) Laplacian. 
By applying our general methods we avoid any specific new work involving the
boundary value problem -- see Theorem \ref{t:3} and Corollary \ref{c:1}.
This and other inverse results are presented in Sect.\ref{s5}.

What may seem like an excessive generality in the trace formul{\ae}
of \cite{SjZw} (see \eqref{eq:trace2} below) and Theorem \ref{t:4} is
motivated by their potential use for the recovery of 
classical dynamics of {\em effective Hamiltonians}. One of
the most impressive inverse procedures is the application of the 
{\em Onsager rules} for determining Fermi surfaces -- 
see for instance \cite{AM}
for an introductory physics discussion. They can be 
interpreted in terms of a trace formula involving an effective Hamiltonian
coming from the {\em Peierls substitution} --  see \cite{dHvA}.
The question which will be investigated elsewhere is: can a more detailed
information about the Fermi surface be obtained from higher order
terms in the expansions of magnetic susceptibilities of metals.

Throughout the paper we will write 
\[ u ( x , \epsilon ) = \sum_{ j = 0}^\infty u_j ( x ) \epsilon ^j \,,\]
to denote asymptotic expansions as $ \epsilon \rightarrow 0 $. When no
confusion is likely to arise all equalities are meant modulo $ {\mathcal 
O} ( \epsilon^\infty ) $ where $ \epsilon $ is the relevant small 
parameter (mostly $ h $, the mathematical ``Planck constant'').
We denote by $ \NN $ the set of natural numbers, $ 0,1,2, \cdots$.

\medskip

\noindent
{\sc Acknowledgements.} The third author would like to thank the
National Science Foundation for partial support under grant DMS-9970614.
He would also like to thank Bill Miller and Howard Taylor for 
references to the chemistry literature.

\section{Preliminaries}
\label{s2}

In this section we will recall the general trace formula of \cite{SjZw}
and the classical and quantum Birkhoff normal forms of \cite{IS}.

We consider $ X $ which is either a compact $ \CI $ manifold of dimension
$ n+1$ or $ \RR^{n+1} $. 

We introduce the usual class of semi-classical 
symbols on $ X $:
\[ S^{m,k} ( T^* X ) = \{ a \in \CI( T^* X \times (0, 1]  ) :
|\partial_x ^{ \alpha } \partial _\xi^\beta a ( x, \xi ;h ) | \leq
C_{ \alpha, \beta} h^{ -m } \langle \xi \rangle^{ k - |\beta| } 
\} \,, \]
and the corresponding class of pseudodifferential operators, $ 
\Psi_{h}^{m,k} ( X ) $, with the quantization and symbol maps:
\[ 
\begin{split}
& \Op \; : \; S^{m , k } ( T^* X ) \ \longrightarrow 
\Psi^{m,k}_h ( X) \\
& \sigma_h \; : \; \Psi_h^{m,k} ( X ) \ \longrightarrow 
S^{ m , k } ( T^* X ) / S^{ m-1, k-1} ( T^* X ) \,, \end{split}
\]
with both maps  surjective,  and the usual properties 
\[ \sigma_h ( A \circ B ) = \sigma_h ( A)\sigma _h ( B ) \,,\]
\[ 0 \rightarrow \Psi^{m-1, k-1} ( X) \hookrightarrow \Psi^{m, k} ( X)
\stackrel{\sigma_h}{\rightarrow} S^{ m , k } ( T^* X ) 
/ S^{ m-1, k-1} ( T^* X ) \rightarrow 0 \,, \]
a short exact sequence, and
\[ \sigma_h \circ \Op : S^{m,k} ( T^* X ) 
\ \longrightarrow  S^{ m , k } ( T^* X ) / S^{ m-1, k-1} ( T^* X ) \,,
\]
the natural projection map.

Let $ P (z ) \in \CI ( I_z ; \Psi_h ^{0, k} ( X ) ) $, 
$ I = ( -a , a ) \subset \RR $, be a family of self-adjoint, 
principal type operators, such that $ \Sigma_z = \{ m \; : \; 
\sigma ( P ( z) ) = 0 \} \subset T^* X $ is compact.
We assume that
\begin{gather}
\label{eq:5.ell}
\begin{gathered}
\sigma ( \partial_ z P ( z) ) 
\leq - C <  0 \,, \ \text{ near $ \Sigma_z $\,, } \ \ \ 
| \sigma ( P ( z) ) | \geq C \langle \xi \rangle^k \,, 
\\ \text{ for $ |\xi | \geq C $ when $ X $ is a compact manifold, and for $ | ( x , \xi ) | \geq C $ if
$ X = \RR^{n+1} $.}
\end{gathered}
\end{gather}
We also assume  that for $ z $ near $ 0$, the Hamilton vector field, 
$ H_{p(z)} $, $ p(z) = \sigma ( P ( z ) )$, has a simple closed
orbit $ \gamma(z) \subset \Sigma_z $ with period $ T ( z) $, 
and that $ \gamma ( z) $ has
a neighbourhood $ \Omega $ such that 
\begin{equation}
\label{eq:t2}
 m \in \Omega \ \text{and} \ \exp t H_{p(z)} ( m ) = m \,, \
p (z,  m ) = 0 \,, \ 0 < |t| \leq T (z) N 
+ \epsilon \,, \ z \in I\,, 
\ \Longrightarrow \ m \in \gamma(z) \,, \end{equation}
where $ T (z) $ is the period of $ \gamma( z) $, assumed to depend smoothly on
$ z $. We also write $ p ( z , m ) $ for $ \sigma ( P ( z ) ) ( m ) $.
Let $ A \in \Psi^{0,0}_h ( 
X ) $ be a microlocal cut-off to a sufficiently small 
neighbourhood of $ \gamma(0) $.

\newcommand{\KKer}{\operatorname{ker}_{m_0 ( z) } ( P ( z ) ) }
\newcommand{\KKert}{\operatorname{ker}_{ \tilde m_0 ( z) 
} ( \widetilde P ( z ) ) }

Then if $ P( z) $ is 
 an almost analytic extension of $ P( z) $, $ z \in \RR $,
$ \chi \in \CIc ( I ) $, $ \tilde \chi \in 
\CIc ( \CC) $, its almost analytic extension, 
$ f \in \CI ( \RR ) $, and $ {\rm{supp}}\; \hat f \subset 
( - N ( C_p - \epsilon)   + C ,  N ( C_p -  \epsilon)  - C )  
\setminus \{ 0 \} $, we have,
\begin{gather}
\begin{gathered}
\label{eq:trace2}
\frac1\pi  {\mathrm{tr}} \;  \int f ( z / h ) 
 \bar \partial_z \left[ \tilde \chi ( z )\;   \partial_z P (z) \;  P ( z) ^{-1} \right] A 
{\mathcal L} ( d z ) 
 = \\
- \frac{1 }{ 2\pi i } 
\sum_{
-1\ne 
k=  -N -1  }^{N -1 }  {\mathrm {tr}} \; 
\int_{\mathbb R} f (z/h) M(z, h)^k \frac{d}{dz} M ( z , h )
 \chi ( z) 
dz + {\mathcal O} ( h^\infty ) \,, \end{gathered} \end{gather}
where $ M ( z, h ) $ is the quantum monodromy operator. 
The constant $  C_p > 0 $, in the condition on $ \hat f $ depends on 
$ p (z) $. 

The {\em quantum monodromy operator}, $ M ( z , h ) $ is
defined as follows.  For a point on $ \gamma $, $ m_0 \in \gamma $, we 
can define the microlocal kernel of $ P ( z ) $ at $ m_0 $, to be the
set of families $ u ( h ) $, such that $ u ( h ) $ are microlocally 
defined near $ m_0 $ and $ P ( z ) u ( h ) = {\mathcal O} ( h^ \infty ) $ 
near $ m_0 $. We denote it by $ \KKer $. Any solution can be 
continued microlocally along $ \gamma ( z) $ and we denote the corresponding
forward and backward continuation by $ I_\pm ( z ) $ -- see 
\cite[\S 4]{SjZw} for precise definitions.  
Now, let $ m_1 \neq m_0 $ be another point 
on 
$ \gamma ( z) $. We then define
\begin{gather}
\label{eq:mon}
\begin{gathered}
 I _{-}  ( z ) {\mathcal M} ( z ) = I_+ ( z ) \,, \ \ \text{ near $ m_1 $,} \\
{\mathcal M} ( z) \; : \;   \KKer
\ \longrightarrow \  \KKer \,. \end{gathered}
\end{gather}

The operator $ P ( z )  $ is assumed to be self-adjoint with respect to some
inner product $ \langle \bullet, \bullet \rangle$, and we define
the {\em quantum flux} norm on $ \KKer $ as 
follows: let $ \chi $ be 
a microlocal cut-off function supported near
$ \gamma $ and equal to one near the part of $\gamma $ between $ m_0 $ and 
$ m_1 $. We denote by $ [ P ( z) , \chi ]_{+} $ the part of the
commutator supported in $ m_0 $, and  put
\[  \langle u , v \rangle_{\rm{QF}}  \stackrel{\rm{def}}{=} 
\langle [ ( h / i ) P ( z) , \chi ]_{+} u , v \rangle \,, \ \ u ,v 
\in \KKer \,. \]
It is easy to check that this norm is independent of the choice of 
$ \chi $  -- see \cite[Lemma 4.4]{SjZw}.  This independence leads
to the unitarity of $ {\mathcal M} ( z ) $:
\[  \langle {\mathcal M} ( z)  u ,  {\mathcal M} ( z )  u \rangle_{\rm{QF}}
=  \langle u , u  \rangle_{\rm{QF}} \,, \ \ u \in \KKer  
\,. \]
For practical reasons we identify $ \KKer $ with 
$ {\mathcal D}' ( \RR^{n} ) $, microlocally near $ ( 0 , 0 )$, and
choose the identification so that the corresponding monodromy map is
unitary (microlocally near $ ( 0 , 0 ) $ where $ ( 0, 0 )$ corresponds
to the closed orbit intersecting a transversal identified with 
$ T^* \RR^{n} $). This gives
\[ M ( z , h ) \; : \; {\mathcal D} ' ( \RR^{n} ) \ \longrightarrow
\ {\mathcal D}' ( \RR^{n} ) \,, \]
microlocally defined near $ ( 0, 0 ) $ and 
\begin{quote} 
$ M ( z , h ) $  is a semi-classical Fourier
integral operator which quantizes the Poincar{\'e} map of $ \gamma ( z ) $.
This, and \eqref{eq:trace2}, \eqref{eq:trace3} are the {\em only} 
properties  of $ M ( z , h ) $ which will be used in this paper.
\end{quote}

Basic ideas become clear when 
one considers the simplest example $ P ( z ) = h D_x - z $, on 
$ {\mathbb S}^1 $. In that case $ M ( z, h ) = \exp ( 2 \pi i z / h ) $
-- see \cite[\S 2]{SjZw} for a careful presentation.

We will now assume that the symbol of $ P ( z ) $ has an asymptotic 
expansion in powers of $ h $, near $ \Sigma_0 $:
\begin{equation}
\label{eq:clas}
(P( z))( x, \xi ; h ) = \sum_{ j = 0 }^\infty h^j (p_j ( z))(x, \xi) \,, \ \
( x , \xi ) \in U\,, \  \ \text{ $ U$ a precompact neighbourhood of 
$ \Sigma_0 $} \,.
\end{equation}

Under this assumption we have a simple extension of 
\cite[Proposition 7.5]{SjZw}
\begin{equation}
\label{eq:trace3}
{\mathrm {tr}} M( z , h )^{k-1} hD_z M( z , h ) =
 e^{ i k I( z) }\frac{ e^{  i \nu_k  ( z) \frac{\pi}{2} }
I' ( z ) 
} {| d \kappa_{ \gamma( z ) } ^k  - 1 |^{\frac12}} 
+  e^{ i k I( z) } \sum_{ j=1}^\infty a_{ j, k } h^{ j } \,,
\end{equation}
where $ I ( z) $ is the classical action, $ \nu_k $ the Maslov index of
the $k$th iterate of $ \gamma ( z) $ (it is independent of $ z $ for
$ z $ near $ 0$).

We start by recalling the general statement 
of the classical {\em Birkhoff normal form}
 -- see \cite[Theorem 1.3, Proposition 4.3]{IS} 
and references there.

Let $ ( W , \omega ) $ be a symplectic manifold of dimension $ 2 n$,
and let
$ \kappa \; : \; W \; \longrightarrow \; W \,, \ \ \kappa ( w_0 ) = w_0 \,,$ 
be a symplectic transformation, $ \kappa^* \omega = \omega $, 
fixing $ w_0 \in W  $.  
In \eqref{eq:ellbnf}, the {\em twist map}, $ \kappa_ 0 $, is the
Birkhoff normal form of $ \kappa $, 
\begin{equation}
\label{eq:bnf1} \kappa_0 = \exp H_p  \,, \
\end{equation}
where the Hamiltonian $ p $ is in the Birkhoff normal form. 

The general case presented in \cite{IS} is formulated using \eqref{eq:bnf1}:
suppose that the eigenvalues of $ d \kappa ( w_0 ) $ are given by 
\begin{gather}
\label{eq:ev}
\begin{gathered}
 \lambda_j \,, \lambda_j^{-1} \,, \
\ j = 1, \cdots, n \\
\lambda_j  = \exp \mu_j \,, \ \ \mu_j = \alpha_j + i \beta_j \,, \ 
\alpha_j >0 \,, \ 0<  \beta_j \notin 2\pi \ZZ   \,, \ \ 
 1 \leq j \leq 
n_{ \rm{ch}} \,, \\
\lambda_{j + n_{\rm{ch}} }  = \exp {\mu_{j+n_{\rm ch}}} \,,\ \mu
_{j+n_{{\rm ch}}}=\overline{\mu _j}\, , \ \ 
1 \leq j \leq 
\NC \,, \\ 
\lambda_j=\exp \mu _j\, ,\ \mu _j> 0 \,, \ \ 
2n_{\rm{ch}} + 1 \leq j \leq 2 n_{\rm{ch}} + n_{ \rm{rh}} 
\,,\\
\lambda_j = \exp \mu_j \,, \ \ 
\mu_j \in i (\RR \setminus  2\pi \ZZ ) \,, \ \ 
n_{ \rm{rh}} + 2 n_{ \rm{ch}} + 1 \leq j \leq n \,.  
\end{gathered}
\end{gather}
Here $ n_{\rm{e}} + n_{\rm{rh}} + 2 n_{ \rm{ch}} = n $, and $ 
\rm{e} $ stands for elliptic, $ \rm{rh} $ for real hyperbolic, and
$ \rm{ch} $, for complex hyperbolic.

We replace the condition in the purely elliptic case \eqref{eq:dio1} by 
the more general condition
\begin{equation}
\label{eq:dio2}
\sum_{ j=1}^n k_j \mu_j \in 2 \pi i \ZZ 
\,, \  
k_j \in \ZZ \ \Longrightarrow \ k_j = 0 \,.
\end{equation}
Then the generalization of \eqref{eq:ellbnf}
says that there exist symplectic coordinates $ ( x, \xi ) \in \RR^n $ 
centered at $ ( 0 , 0 ) $ 
\begin{gather}
\label{eq:bnf}
\begin{gathered}
\kappa = \exp H_p + \kappa^\flat  \,, \ \ \kappa^\flat
 ( x , \xi ) = {\mathcal O} ( ( x^2 + \xi^2 )
^\infty ) \,, \\ 
p = p_0 + r \,, \ \ r ( x , \xi ) = {\mathcal O} ( ( x^2 + \xi^2) ^2 ) \,, \\
 \begin{split} 
& p_0 ( x, \xi ) =  
  \sum_{ j = 1 }^{ \NC} \alpha_j ( x_{ 2j - 1} \xi_{ 2j - 1}
+ x_{2j } \xi_{ 2 j } )  -  
\beta_j  ( x_{ 2j - 1} \xi_{ 2j }
-  x_{2j } \xi_{ 2 j -1 } ) \\
& + \; 
\sum_{ j = 2 \NC + 1}^{ 2 \NC + \NR } \mu_j x_j \xi _j  \ +  
\sum_{ j = 2 \NC + \NR + 1 }^{n} \frac{\mu_j }{ 2 i } ( x ^2_j + \xi_j^2 )  
 \end{split}
  \\
r ( x , \xi ) = R ( \iota ) \,, \ \ R ( \iota ) = \sum_{ |\alpha | \geq 2 }
r_\alpha \iota^\alpha \,, \\
 \iota_j =  \left\{ 
\begin{array}{ ll} 
(  x_{ 2j - 1} \xi_{ 2j - 1} + x_{2j } \xi_{ 2 j } 
+ i   (x_{ 2j - 1} \xi_{ 2j } -  x_{2j } \xi_{ 2 j -1 })) /2  &
 1 \leq j \leq  \NC \\
\bar \iota_{j'}  & 
j' = j - \NC \,,
\   \NC + 1 \leq j \leq 2 \NC  \\
x_j \xi_j & 2 \NC + 1 \leq j \leq 2 \NC + \NR \\
 ( x ^2_j + \xi_j^2)/(2i)   & 2 \NC + \NR + 1 \leq j \leq n \,. 
\end{array} \right.
\end{gathered}
\end{gather}
Observe that $ p_0 = \langle \iota, \mu \rangle $.
By the {\em Birkhoff normal form} of $ \kappa $ we mean $ p_0 $ and
the asymptotic expansion of $ r $.

We now recall the {\em quantum} Birkhoff normal forms \cite{I} and
\cite[Theorem 5.1, (4.11)]{IS}. We start with some general facts about
$ h$-Fourier integral operators. Let $ W $ above be a neighbourhood of
$ w_0 \in T^* X $, where $ X $ is a $ \CI $ manifold. To $ \kappa : W 
\rightarrow W $ we can associate a class of operators, microlocally 
defined near $ w_0 $ (see \cite[\S 3]{SjZw} for the precise definition of
this notion):
\[  I_h^m ( X , X , \kappa' ) \,,\]
the semi-classical {\em Fourier integral operators} quantizing $ \kappa $.
One way to define them is by considering oscillatory integrals: we identify 
$ w_0 $ with $ ( 0 , 0 ) \in T^* \RR^n $ so that
\begin{gather}
\label{eq:fio}
\begin{gathered}
U \in I_h^m ( X , X , \kappa' )  \ \Longleftrightarrow \ 
U u ( x ) = h^{ -m - \frac{n+N} 2 } \int \int  e^{ i \phi ( x, y , \theta ) /h }
a ( x , y , \theta ; h ) u ( y ) d y d \theta \,, \\
\phi \in \CI ( \RR^n \times \RR^n \times \RR^N ) \,, \ \ 
{\rm{graph}}( \kappa ) = \pi ( C_\phi ) \,, \\
C_\phi = \{ ( x, y , \theta ) \; : \; \phi_\theta' ( x , y , \theta ) 
= 0 \} \ni ( x, y , \theta ) \stackrel{\pi}{\longmapsto} ( x , \phi'_x ;
y , \phi_y' ) \,, \end{gathered}\end{gather}
where $ \phi $ is assumed to be a non-degenerate phase function in the sense of
H{\"o}rmander (see \cite[Def. 21.2.5]{Horb} for a related discussion), 
that is $ \phi $ is smooth and real, $ C_\phi $ is a smooth manifold of
dimension $ 2n $, and $ \pi $ has an injective differential.
The amplitude $ a $ is assumed to be supported near $ ( 0 , 0, 0 ) $, and
to be a classical symbol in the sense of \eqref{eq:clas}.

Before stating the main result of \cite{IS} we have to introduce the
notion of equivalence of families of 
operators used in that paper \cite[\S 2, \S 5]{IS}:
let $ U ( z) \in I^0 ( X , X , \kappa_z' ) $ and $ \widetilde U ( z) \in 
I^0 ( X , X , \tilde \kappa_z' ) $ be two families of operators. 
Then 
\begin{equation}
\label{eq:eq} U \equiv \widetilde U 
\end{equation}
 if 
$ \kappa_0 ( w_0 ) = \tilde \kappa_0 ( w_0 ) = w_0 $, the two families
of transformations, $ \kappa_z $ and $ \tilde \kappa_z $,
agree to infinite order at $ w_0 $ and $ z = 0 $, and the terms in 
the asymptotic expansions (in powers of $ h$) of $ \phi_z $, 
$ \tilde \phi_z $, $ a_z $, $ \tilde a_z $
(with amplitudes as in \eqref{eq:fio}) agree to infinite order at 
$ ( 0 , 0 , 0 )$ and $  z = 0 $. 

Suppose that $ U ( z) $ quantizes $ \kappa_z  $ 
and $ \kappa_0 $  satisfies the assumptions
of \eqref{eq:bnf}. Also assume that $U(z)$ is elliptic. Then
\begin{gather}
\label{eq:qbnf}
\begin{gathered}
U ( z) \equiv V( z) ^{-1} \, e^{ i P ( z , h ) / h } \, V( z)  \,, \\ 
V ( z) \in I^0 ( X ,  X, 
{\mathcal C}_z' )  \ \text{ is a family of
 elliptic Fourier integral operators near $ ( (0, 0 ) , w_0) $}
 \\
P ( z, h ) = p_0 ( z) + hp_1(z)  + h^2 p_2 ( z) 
+ \cdots \,, \\ p_j ( z)  = p_j (z;  I_1, \cdots, I_j )
\,, \ \ I_j =  \iota^w_j ( x , h D_{x} ) \,, \\ 
p_0  ( z) = \sum_{ j = 1}^{ n } \mu_j ( z) I_j  + R ( z; 
I_1 , \cdots , I_n ) \,, \ \
R ( z; \iota ) = {\mathcal O} ( \iota^2 )  \,. \\
\end{gathered}
\end{gather}
When $ U ( z) $ is microlocally unitary near $ w_0 $, $P( z , h) 
 $, $ V( z)  $ can be chosen
to be microlocally self-adjoint and unitary, respectively, 
at $ ( (0, 0 ) , w_0)  $.
We will call the expansion of $ P ( z , h ) $ at $ z=0 $, $ \iota = 0 $,
the {\em infinitesimal Quantum Birkhoff Normal Form} of $ U ( z )$ at $ z = 0 $.

We recall that a Fourier integral operator is elliptic if it is 
associated to a canonical transformation, that is $ {\mathcal C} $ above
is the graph of a canonical transformation, $ 
( (0, 0 ) ,  w_0  ) \in {\mathcal C} $, and that its symbol
is elliptic, that is the leading term in the
asymptotic expansion of the amplitude in  \eqref{eq:fio} does not vanish.

\section{Recovering Birkhoff normal forms from traces}
\label{s3}

From \cite[Lemma 3.2, Proposition 7.3]{SjZw} we know that if 
$ U (z) \in I^0_h ( X , X , \kappa( z)' ) $, is a family of 
$h$-Fourier integral operators,  and that the canonical relations
$ \kappa ( z)^k $, $ k \ne 0$,  have only one
non-degenerate fixed point then 
\begin{equation}
\label{eq:tr}
\tr  \; U ( z) ^k  = e^{ i k S( z)  /h } \left( a_{0k}( z)  + a_{1k} ( z) 
h + \cdots + a_{lk} ( z) 
{ h^l} + {\mathcal O} ( h^{ l+1} ) \right) \,.\end{equation}
In this section we will prove the following

\medskip

\begin{thm}
\label{t:1}
Let $ U ( z) \in I^0_h (
X , X , \kappa(z)' ) $ be a family of
 $ h$-Fourier integral operator associated to 
smooth locally defined 
canonical transformations, $ \kappa ( z) : T^*X \rightarrow T^*X $.
Suppose that $ \kappa_0 $ has a unique fixed point, $ w_0$, and that
the eigenvalues of 
$ d \kappa_0 ( w_0 ) $ given in \eqref{eq:ev} satisfy \eqref{eq:dio2}. 
Then 
the Taylor series of the coefficients $ a_{lk} ( z) $ in \eqref{eq:tr} 
determine the  {\em infinitesimal
quantum Birkhoff normal form} of $ U ( z) $, $ P ( z , h )$ in \eqref{eq:qbnf}.
In particular, the 
Birkhoff normal form of $ \kappa_0 $
at $ z = 0 $, $ p ( z) $, in \eqref{eq:bnf}, is determined.
\end{thm}

\medskip

Our proof of Theorem \ref{t:1} is based on 
\begin{prop}
\label{p:3.1}
Suppose that 
\begin{gather}
\label{eq:exp1}
\begin{gathered}
G ( \iota ; h ) =  \langle \iota , \mu \rangle + F ( \iota, h )\,, \ \ 
F ( \iota, h ) =  F_0 ( \iota ) + h F_1 ( \iota ) 
+ \cdots \,, 
 \ \ e^{ \mu_j } \neq 1 \,, \ \
F_0 ( \iota ) = {\mathcal O} ( \iota^2 ) \,, 
\end{gathered}
\end{gather}
where we consider $ F ( \iota, h ) $ as a formal power series in $ h$, and
each $ F_j ( \iota ) $, a formal power series in $ \iota $, and where
we supressed the dependence on $ \mu $ in higher order terms.

If $ I_j ( h ) = \iota_j ^w ( x , h D_x ) $, then 
\begin{equation}
\label{eq:tr1}
\tr \; e^{ -  i G ( I ( h ) ; h ) /h} =
e^{ -i  F ( i h \partial_\mu ; h ) / h } \left( \prod_{ j=1}^{n} 
\frac{1}{ 2 \sinh (  \mu_j /2 ) } \right) \,.
\end{equation}
\end{prop}
\begin{proof}
When $ G ( \iota ; h ) = \langle \iota , \mu \rangle $ this follows from 
computations in each of the cases given in \eqref{eq:bnf}:
\[ \tr \; e^{- i \langle I ( h ) , \mu \rangle / h } = 
\tr \; e^{ -i \langle I ( 1) , \mu \rangle } \,, \]
where we can separate the trace into separate blocks:
\begin{itemize}

\item {\em elliptic blocks:} $$  \tr e^{ - \mu ( D_x^2 + x^2 )/2 } = 
\sum _{ k = 0}^\infty e^{ - \mu ( n + \frac12 )} = \frac1{ 2 \sinh \frac{\mu} 2 
} \,, $$
where the formal computation ($ \mu \in i \RR_+ $) can be justified by 
condsidering the limit $ \epsilon \rightarrow 0 $, with $ \mu $ 
replaced by $ \mu + \epsilon$.

\item {\em real hyperbolic blocks:} $$ \tr e^{ - i \mu ( x D_x + D_x x ) /2 }
= \tr e^{ - \mu ( x \partial_x + \frac12 ) } = \tr \left( 
u \mapsto e^{ -\mu/2 } u ( e^{ -\mu }\bullet ) \right) \,, $$
Formally, 
\[ \tr e^{ - \mu ( x \partial_x + \frac12 ) } = e^{ - \mu/2 } 
\int \delta_0 ( ( e^{-\mu } -1 ) x ) dx = \frac{1}{ 2 \sinh \frac{\mu}2} \,,
\]
and this can be justified by approximating $ \delta_0$ by a sequence of 
smooth compactly supported functions.

\item {\em complex hyperbolic blocks:} here we have to work in two dimensions,
$ \mu_1 = \alpha + i \beta $, $ \mu_2 = \bar \mu_1 = \alpha - i \beta $, 
and the corresponding $ \iota_1 $ and $ \iota_2 $ (see the last two lines in 
\eqref{eq:bnf}) are $ \iota_1 = z \cdot \zeta $, $ \iota_2 = w \cdot \omega $,
where $ z = ( x_1 - i x_2 )/\sqrt 2 $, $ w= \bar z $ and 
$ \zeta =
( \xi_1 + i \xi_2) / \sqrt 2 $, $ \omega = ( \xi_1 - i \xi_2 )/ \sqrt 2 $, 
the corresponding dual coordinates, $ \omega = \bar 
\zeta $. With $ I_j = \iota_j^w ( x, D_x ) $, $ j = 1, 2$, we have
\[ \tr e^{ - i \mu_1 I_1 - i \mu_2 I_2 } = 
\tr e^{ - \mu_1 ( z \partial_z + \frac12 )  - \bar \mu_1 ( \bar z \partial _{
\bar z } + \frac12 ) }  = \tr \left( u ( \bullet , \bar \bullet ) 
\mapsto e^{ - \mu_1/2 - 
\bar \mu_1/2 } u ( e^{ - \mu_1 } \bullet , e^{ - \bar \mu_1 } \bar \bullet )
 \right) \,.\]
The last operator can be written as a contour integral,
\[  e^{ - \mu_1/2 - 
\bar \mu_1/2 } u ( e^{ - \mu_1 }  , e^{ - \bar \mu_1 } \bar z  ) =
\frac{1}{ ( 2 \pi)^2 } \int \! \!  \! \int \! \! \! \int \! \!\! \int e^{ i ( ( e^{- \mu_1 } z - \tilde z )
\tilde \zeta + ( e^{- \bar \mu_1 } \bar z -  {\tilde w } ) \tilde \omega )}
u ( \tilde z , \tilde w ) d \tilde z d \tilde w d \tilde \zeta d \tilde 
\omega  \,,\]
with the integration contour given by 
$  \tilde w = \bar{ \tilde z } \,, \ \tilde \omega = \bar{ \tilde \zeta } $.

Formal stationary phase argument applied when taking the trace gives
$$ \prod_{ j =1}^2 \frac{1}{ 2 \sinh \frac{ \mu_j}2 } \,. $$
\end{itemize}

Once we understand the contribution of the leading term 
the following  formal argument essentially gives \eqref{eq:tr1}. Since
$ I ( h ) $ is unitarily equivalent to $ h I ( 1 ) $ we see that
\[ \tr \; e^{ -  i G ( I ( h ) ; h ) /h} = 
\tr \; e^{ -  i G ( h I ( 1 ) ; h ) /h} \,. \]
On the other hand, $ I( 1)  = \partial_\mu  \langle I( 1) , \mu  
\rangle $, and hence 
\[ e^{ - i G ( h I ( 1 ) ; h ) / h } = 
 e^{ -i F ( h I ( 1 ) ; h )/ h} e^{ -i \langle  I ( 1 ) , \mu \rangle  } =
e^{ -i F   ( i h \partial_\mu ; h ) / h } e^{ -i \langle I ( 1 ) , \mu  
\rangle  }
\,. \]
By moving differentiation outside of the trace we obtain the expansion
\eqref{eq:tr1}.
\end{proof}

We also need a result of Fried \cite{F} used in the same context in 
\cite{G}:
\begin{prop}
\label{p:3.2}
If the assumption \eqref{eq:dio2} holds then the set
\[ a_0^k  \prod_{ j=1}^{n}   \sinh ( k \mu_j /2 )  \,, \ \ k \in \NN \,,
\]
determines $ \exp \mu_j$'s and $ a_0 $ uniquely.
\end{prop}
\stopthm

Proposition \ref{p:3.2} and an argument based on Kronecker's 
Theorem give the following inverse result crucial in the proof of 
Theorem \ref{t:1}  (see  \cite[\S 8]{G}, \cite[\S 6]{Zel2}) for different 
original approaches):

\begin{lem}
\label{l:3.1}
Suppose that $ \mu_j $ satisfy \eqref{eq:dio2}.
Then for any 
polynomial,  $ p ( \xi) =
 \sum_{|\alpha| \leq M } a_\alpha \xi^\alpha $, the coefficients
$ a_\alpha $ can be recovered from the asymptotics of
\begin{equation}
\label{eq:l3}
 p ( i k^{-1} \partial_\mu ) \prod_{ j =1 }^{n} \frac{1}{ 2 \sinh ( k \mu_j /2 ) }  
\,, \  \  \ k \longrightarrow \infty \,. \end{equation}
\end{lem}
\begin{proof}

Using the argument of the proof of Proposition \ref{p:3.1} we can rewrite
\eqref{eq:l3} with $e_0=(1,...,1)$: 
\begin{equation}
\label{eq:l3'}
\begin{split}
    p ( i k^{-1}  \partial_\mu ) \prod_{ j =1 }^{n} \frac{1}{ 2 \sinh 
( k \mu_j /2 ) }  
& =      p ( i k^{-1} \partial_\mu ) 
\sum_{ m\in \NN^{n}  }  e^{ - \langle m + e_0/2 , k \mu 
\rangle } \\
& = \sum_{ m\in \NN^{n}  } p (( m + e_0/2)/i ) e^{ - \langle m + e_0/2 , k \mu 
\rangle } \,. 
\end{split}
\end{equation}
To justify this expression we will consider it as a distribution 
obtained as a boundary value of a holomorphic function. We introduce
new complex variables $ z \in \Omega_1 = D(0,1)^{n_{\rm{e}}} $ and $ 
w \in \Omega_2 =( D( 0 , 1) \setminus {0} )^{n_{\rm{ch}} } $ in a polydisc and
a punctured polydisc respectively: at the boundary, 
$ z_j = e^{  i \mu_j/2 } $, and 
$ w_j  = e^{ i  \beta_{j + n_{\rm{e}} + n_{{\rm{rh}} }}/2 } $. 
We then write 
\begin{gather*}
m = ( m_{\CH}' , m_{\CH}'' , m_{\RH} , m_{\EE} ) \in \NN^n \,, \ \ 
e_0 = ( e_0 ^{\CH}, e_0^{\CH} , e_0^{\RH} , e_0^{\EE} ) \,, \\
\mu = ( \mu_{\rm{ch}}, \bar \mu_{\rm{ch}},  \mu_{ {\RH}}, \mu_{\rm{e}}) \,,
\ \  \mu_{\rm{ch}} = \alpha + i \beta \,, \ \ 
\mu_{\RH} \in \RR_+^{n_{\RH}} \,, \ \
\mu_{\EE} \in ( i \RR)^{n_{\EE}} \,, \\
\langle m + e_0/2 , \mu \rangle = 
\langle n_2 +e_0^{\CH} , \alpha \rangle + \langle n_1 + e_0^{\RH} /2 , 
\mu_{\RH} \rangle + i \langle l_1 + e^{\EE} / 2 ,  \mu_{\EE}/ i \rangle +
i \langle l_2 , \beta \rangle \,,\\
n_2 = m'_{\CH} + m''_{\CH} \,, \ \ l_2 = m_{\CH}' - m_{\CH}'' \,, \ \ 
n_1 = m_{\RH } \,, \ \ l_1 = m_{\EE} \,. 
\end{gather*}
Here $ m \in \NN^n $ is arbitrary, and
\[ n_1 \in \NN^{\NR} \,, \ \ l_1 \in \NN^{\NE} \,, \ \
n_2 \in \NN^{\NC} \,, \ \ l_2 \in \ZZ^{\NC} \,, \]
with the constraint that 
\[  \frac{ n_2 \pm l_2 }{2} \in \NN^{\NC} \,, \]
which implies that $ | l_2 | \leq n_2 $, component wise. We then put
\[ q ( n , l ) = p ( ( m + e_0 /2 ) / i )\,, 
\ \ z = e^{ - \mu_{\EE} } \in \CC^{ {\NE} } \,, \ \ w = e^{ - i \beta
  / 2}\in \CC^{\NC} \,,
\]
(component-wise in the natural sense). In this notation \eqref{eq:l3'}
becomes
\begin{equation}
\label{eq:l3''}
 \sum_{ n \in \NN^{n_{\rm{rh}} + n_{\rm{ch}} } }
e^{ - k (\langle n_1 + e^{\rm{rh}}_0/2 , \mu_{\rm{rh}}  \rangle + 
\langle n_2 + e^{\rm{ch}}_0/2 , \alpha \rangle )} 
F_{ n } ( z^k , w^k ) \,,  \end{equation}
where 
\[ 
F_n ( z , w ) = \sum_{{ l_1 \in  \NN^{n_{\rm{e}} } \,, \ 
 l_2 \in \ZZ^{n_{\rm{ch}} } }\atop{   \frac{n_2 \pm l_2 }2 \in \NN^{\NC} }  }
q ( n , l ) z^{ 2 l_1 + e^{\rm{e}}_0 }  
w^{ 2 l_2  } \,, \]
$  ( z , w ) \in \Omega_1 \times \Omega_2 \subset \CC^{n_{\rm{e}}} \times
\CC^{n_{\rm{ch}}} $. Here we note that for any fixed $ n $, $ q ( n , l ) $
is non-zero for finitely many $ l_2$'s and hence $ F_n( z , w )$ 
is holomorphic in the 
punctured polydisc. Our task to recover the coefficients $ q ( n , l ) $ from 
asymptotics as $ k \rightarrow 0 $. That amounts to recovering the 
holomorphic functions $ F_n $. We observe that \eqref{eq:l3'} shows that
the distributional boundary values  $ F_n $'s are smooth away from $ z_j^4 =1 $.

An application
of Kronecker's theorem\footnote{which says that if $ \alpha \in \ZZ^n $
satisfies $ \langle \alpha,  k\rangle \in 2 \pi \ZZ $, $ k \in \ZZ^n $, 
$ \Rightarrow \ k = 0 $, 
then for any $ y \in  \RR^n $ and $ \epsilon > 0 $,
there exists $ N \in \NN $ and $ K \in \ZZ^{ n } $, such that
$ | y - N \alpha  - 2 \pi   K | < \epsilon.$ Its use is
justified by \eqref{eq:dio2}.} 
shows that for any $$ ( x , y ) \in 
(\partial D(0,1))^{n_{\rm{e}}}  \times (\partial D(0,1))^{n_{\rm{ch}}} \,,$$ 
there exist a sequence $ k_r \rightarrow \infty$ such that 
$ ( z^{k_r} , w^{k_r} ) \rightarrow ( x , y ) $.

From  \eqref{eq:dio2} we deduce that
\begin{equation*}
\sum_{ i =  1}^{ 
 n_{ \rm{rh}}  + n_{\rm{ch}} } 
 k_j \Re \mu_{j } = 0 \ \Longrightarrow k_j = 0 \,.
\end{equation*}
Hence by choosing 
sequences of $ k \rightarrow \infty $,
we can determine the coefficients, $ 
F_n ( x, y ) $ in the first sum \eqref{eq:l3''} (the terms have different
rates of exponential decay). From that we determine the holomorphic functions
$ F_n $, and the coefficients, $ q ( n , l ) $.
\end{proof}

\medskip

\noindent
{\em Proof of Theorem \ref{t:1}.} 
Applying \cite[Theorem 3.2]{IS} and \cite[Theorem 4.4]{IS} we see
that, up to $ {\mathcal O}(h^\infty ) $  errors, we can write 
\[ \tr \; U (z ) ^k = \tr \;  e^{ - k i G ( I ( h ), z ; h ) /h} \,, \]
where $ G $ is as in Proposition \ref{p:3.1}. Hence,
\[ \tr \; U ( z ) ^k = e^{ - i k F ( ih k^{-1} \partial_\mu, z  ; h )/h} \tr \; 
e^{ - i \langle I ( 1) , k \mu \rangle } \big|_{\mu = \mu ( z) } 
\,.\]

Recovering the Birkhoff normal form of 
$ \kappa $ means recovering the Taylor series of
$ \langle \iota ,\mu ( z )  \rangle + 
F_0 ( \iota , z )$ from the
trace. We will actually recover the full infinitesimal
{\em quantum} Birkhoff
normal form, $ \langle \iota , \mu ( z) \rangle + F ( \iota,  z; h) $.

We write
\begin{gather*}
 \frac{ F ( h y, z  ; h ) }{ h } = \sum_{ j=0}^\infty h^{j} f_j (z, 
 y ) \,, \ \ 
f_j (z,  y ) = \sum_{ l=0}^{j+ 1} \sum_{ |\alpha | = j +1 - l } 
\frac{1}{ \alpha!} \partial^\alpha F_l ( z, 0 ) y^\alpha  \,, \\
f_j ( z , y ) = \sum_{k=0}^\infty f_{jk} ( y ) z^k  \,.
\end{gather*}
The polynomials $ f_{jk} (  y ) $ clearly determine the formal 
series $ F ( \iota , z ; h ) $.

Since 
\begin{equation}
\label{eq:3.n}
 \begin{split}
& \exp ( -i k F ( h i k ^{-1} \partial_\mu ,  z, ; h ) / h ) = \exp
(-ik\sum_0^{\infty }h^{j}f_j(z,ik^{-1}\partial _\mu ))= \\
& \ \ 
e^{-i k f_0 ( z )  } \left( 1 -  i \sum_{ j=1}^\infty h^{j} k f_j (z,  ik^{-1}
 \partial_\mu  ) -  
\frac{1}{2!} \left( 
 \sum_{ j=1}^\infty h^{j} k f_j (z,  i k^{-1} 
\partial_\mu  ) \right)^2 + \cdots \right)
\,,
\end{split} \end{equation}
we see from \eqref{eq:tr1}, that the coefficient of $ h^0 $ in the
expansion of the trace is given by
\[ \begin{split}
& e^{ -ik f_{00} } \left( 
( 1 -  i \sum_{ m=1}^\infty z^{m} k f_{0m}  -  
\frac{1}{2!} \left( 
 \sum_{ m=1}^\infty z^{m} k f_{0m} \right)^2 + \cdots \right) 
 \   \ \times \\ & 
 \left(  
1 -  i \sum_{ j=1}^\infty z^{j} k \langle  ik^{-1}
 \partial_\mu , \mu^{(j)} ( 0)  \rangle -  
\frac{1}{2!} \left( 
 \sum_{ j=1}^\infty z^{j} k \langle  i k^{-1} 
\partial_\mu  , \mu^{(j)} (0) \rangle  \right)^2 + \cdots \right) 
\tr e^{ -i\langle I ( 1 ) , k\mu \rangle }
 \big|_{\mu = \mu( 0) }  \,, \end{split} \]
where we expanded $ \mu ( z ) $ into its Taylor series and then
used the same argument as in the proof of Proposition \ref{p:3.1}.
We apply Proposition \ref{p:3.2} to the constant term to recover 
$ \exp{i f_{00}}  $ and $ \mu(0)$. 
For $ j > 0 $, the coefficient of $ z^j $ in the expression above 
is given by
\[ -ik e^{ -i k f_0 } (  f_{0j} + \langle k^{-1} i \partial_\mu, 
\mu^{(j)} ( 0 ) \rangle + 
H_{j , k } ( i k^{-1} \partial_\mu) ) 
\left( \prod_{ j =1}^{n} \frac{1}{ 2 \sinh ( k \mu_j  /2 ) } \right)
 \big|_{\mu = \mu( 0) }  \,,
\] 
where $ H_{ j, k } $ is constructed from $ f_{0l} $'s and 
$  \mu^{(p)} ( 0 ) $ with $ l,p < m $. Hence Lemma \ref{l:3.1} shows
that $ f_{0j}$'s and $ \mu^{(j)} ( 0 ) $ can be successively recovered.

We now move to higher powers of $ h$, with the procedure being
essentially the same:
we see from \eqref{eq:3.n} that the coefficient of $ h^m $ in the
expansion of the trace is given by
\[ e^{ -i k f_0 ( z) } (-i  k f_m ( i k^{-1} \partial_\mu , z ) + 
G_{m , k } ( i k^{-1} \partial_\mu, z ) ) 
\left( \prod_{ j =1}^{n} \frac{1}{ 2 \sinh ( k \mu_j /2 ) } \right) 
\big| _{ \mu = \mu ( z) } \,,
\] 
where $ G_{ m, k } $ is constructed from $ f_l $'s with $ l < m $, and
we already know the $ z$-expansions of $ f_0 ( z) $ and $ \mu ( z) $. 
Replacing $ f_m $ by its expansion in $ z $, reduces the problem to 
recovering the polynomials $ f _{ml } ( y ) $ from 
\[  f_{m l}  ( i k^{-1} \partial_\mu )  
\prod_{ j = 1 }^{n} \frac{1}{ 2 \sinh ( k \mu_j /2 ) }  \,, 
\ \ k \in \NN \,, \]
and that follows again from Lemma \ref{l:3.1}. 
This shows that we recovered the full expansion (in $ \iota$, $ z $, and 
$ h $) of $ G ( \iota, z ; h )$, that is the infinitesimal quantum Birkhoff
normal form at $ z  = 0 $.
\stopthm

\section{The monodromy operator determines $P$ near $\gamma $.}
\label{s4}

\renewcommand{\KKer}{\operatorname{ker}_{m_0  } ( P - z  ) }
\renewcommand{\KKert}{\operatorname{ker}_{ \tilde m_0  
} ( \widetilde P -z  ) }

The trace formula \eqref{eq:trace3} and the results of Sect.\ref{s3} will be 
used to recover the infinitesimal quantum Birkhoff normal form 
of the quantum monodromy operator $ M ( z ) $ at $ z = 0$. 
In this section we will discuss the recovery of the full Hamiltonian $ P ( z)$,
infinitesimally at $ z= 0$, from $ M ( z) $. 
Since our information comes from the left hand side of the
trace formula \eqref{eq:trace3}, this recovery can only be possible up to
conjugation by families of elliptic $h$-Fourier integral operators, and
up to multiplication by families of elliptic $h$-pseudo-differential operators.
The assumptions \eqref{eq:5.ell}, the implicit function theorem, and
the usual symbolic iteration, imply that 
\[ P ( z) = A ( z)^*  ( P - z ) A ( z)\,,\]
with $  A ( z) \in \Psi_h ^{0, k/2 } ( X ) $ elliptic near $ \gamma (0)$, and
$ P \in \ \Psi_h ^{0, 0 } ( X ) $ is self-adjoint. Replacing
$ P ( z) $ by $ P - z $ in \eqref{eq:trace3} changes the trace by $ {\mathcal
O }( h^\infty) $. Hence we would like to recover $ P $ infinitesimally 
near $ \gamma $,  on the energy surface $ p ^{-1} ( 0 ) $, 
$ p = \sigma ( P ) $.

\def\l{\ell}
\def\parno{\par \noindent}
\def\vvekv#1#2#3{$$\leqalignno{&{#2}&({#1})\cr
&{#3}\cr}$$}
\def\vvvekv#1#2#3#4{$$\leqalignno{&{#2}&({#1})\cr
&{#3}\cr &{#4}\cr}$$}
\def\vvvvekv#1#2#3#4#5{$$\leqalignno{&{#2}&({#1})\cr
&{#3}\cr &{#4}\cr &{#5}\cr}$$}
\def\ekv#1#2{$${#2}\eqno(#1)$$}
\def\eekv#1#2#3{$$\eqalignno{&{#2}&({#1})\cr &{#3}\cr}$$}
\def\eeekv#1#2#3#4{$$\eqalignno{&{#2}&({#1})\cr &{#3}\cr &{#4}\cr}$$}
\def\eeeekv#1#2#3#4#5{$$\eqalignno{&{#2}&({#1})\cr &{#3}\cr &{#4}\cr
&{#5}\cr}$$}
\def\eeeeekv#1#2#3#4#5#6{$$\eqalignno{&{#2}&({#1})\cr
&{#3}\cr &{#4}\cr &{#5}\cr &{#6}\cr}$$}
\def\eeeeeekv#1#2#3#4#5#6#7{$$\eqalignno{&{#2}&({#1})\cr &{#3}\cr &{#4}\cr
&{#5}\cr&{#6}\cr&{#7}\cr}$$}
\def\cint{{1\over 2\pi i}\int}
\def\iint{\int\hskip -2mm\int}
\def\iiint{\int\hskip -2mm\int\hskip -2mm\int}
\def\buildover#1#2{\buildrel#1\over#2}
\font \mittel=cmbx10 scaled \magstep1
\font \gross=cmbx10 scaled \magstep2
\font \klo=cmsl8
\font\liten=cmr10 at 8pt
\font\stor=cmr10 at 12pt
\font\Stor=cmbx10 at 14pt
\def\aby{arbitrary}
\def\ably{arbitrarily}
\def\asy{asymptotic}
\def\bdd{bounded}
\def\bdy{boundary}
\def\coef{coefficient}
\def\coeff{coefficient}
\def\const{constant}
\def\Const{Constant}
\def\canform{canonical transformation}
\def\coef{coefficient}
\def\coeff{coefficient}
\def\cont{continous}
\def\cor{correspond}
\def\diff{diffeomorphism}
\def\diffeo{diffeomorphism}
\def\de{differential equation}
\def\dop{differential operator}
\def\ef{eigen-function}
\def\ev{eigen-value}
\def\e{equation}
\def\fy{family}
\def\fu{function}
\def\F{Fourier}
\def\fop{Fourier integral operator}
\def\fourior{Fourier integral operator}
\def\fouriors{Fourier integral operators }
\def\hol{holomorphic}
\def\indep{independent}
\def\lhs{left hand side}
\def\mfld{manifold}
\def\ml{microlocal}
\def\mly{microlocally}
\def\neigh{neighbourhood}
\def\nondeg{non-degenerate}
\def\op{operator}
\def\og{orthogonal}
\def\pb{problem}
\def\Pb{Problem}
\def\pde{partial differential equation}
\def\pe{periodic}
\def\pert{perturbation}
\def\Prop{Proposition}
\def\pol{polynomial}
\def\pop{pseudodifferential operator}
\def\pseudor{pseudodifferential operator}
\def\res{resonance}
\def\rhs{right hand side}
\def\sa{self-adjoint}
\def\schr{Schr{\"o}dinger operator}
\def\sop{Schr{\"o}dinger operator}
\def\st{strictly}
\def\stpsh{\st{} plurisubharmonic}
\def\strans{^\sigma \hskip -2pt}
\def\suf{sufficient}
\def\sufly{sufficiently}
\def\sm{\setminus}
\def\tf{transformation}
\def\Th{Theorem}
\def\th{theorem}
\def\tf{transform}
\def\traj{trajectory}
\def\trajs{trajectories}
\def\trans{^t\hskip -2pt}
\def\top{Toeplitz operator}
\def\uf{uniform}
\def\ufly{uniformly}
\def\vf{vector field}
\def\wrt{with respect to}
\def\Op{{\rm Op\,}}
\def\Re{{\rm Re\,}}
\def\Im{{\rm Im\,}}

Thus in this section, 
let $P$ be a \sa{} $h$-\pop{} with leading symbol $p$. Let $\gamma $
be a simple closed $H_p$-\traj{} of period $T_0>0$. Recall from 
\cite[\S 4]{SjZw}
that the quantum time, $ Q$, 
is an $h$-\pop{} microlocally defined in a suitable
\neigh{}
$W_+={\rm neigh\,}(m_0)$ of some fixed point $m_0\in\gamma $, such that
\[
{ {i\over h} [P,Q]}  = Id \,. \] 
After conjugation by a unitary \fop{}, we may assume that
\begin{equation}\label{n:2}
{P=hD_{x_n},\ Q=x_n.}\end{equation}
Let $K(z):{\mathcal D}'({\bf R}^{n})\to {\mathcal D}'(X)$ be a microlocally
defined \fop{} with
\begin{equation}
\label{n:3}
{(P-z)K(z)=0\,, \ \ \  {i\over h} K(z)^*[P,\chi ]_{W_+}K(z)= Id \,, 
}\end{equation}
where $\chi \in C^\infty (W_+)$ is equal to $0$ near the incoming part of
$\gamma $ and equal to 1 near the outgoing part, and we let $\chi $ also
denote a \cor{}ing $h$-quantization. In the case of \eqref{n:2}, we can take
\begin{equation}
\label{n:4}
{K(z)u(x)=e^{izx_n/h}u(x').}\end{equation}
It follows from the second part of \eqref{n:3}, that
\[ {K(z)^*f'(Q)K(z)= Id \,, } \]
for $f\in C^\infty ({\bf R})$, $=0$ for large negative $t$ and $=1$ for
large positive $t$. Indeed, $e^{itP/h}Qe^{-itP/h}=Q+t$ and
\[ \begin{split}
{i\over h} [P,f(Q)]&=(\partial _t)_{t=0}(e^{itP/h}f(Q)e^{-itP/h})=(\partial
_t)_{t=0} f(e^{itP/h}Qe^{-itP/h})\\
&=(\partial _t)_{t=0}f(Q+t)=f'(Q).\end{split} \]
 Notice that $f'\in C_0^\infty ({\bf R})$, $\int f'(q)dq=1$. More generally,
\begin{equation}
\label{n:6}
{K(z)^*g(Q)K(z)=\int g(q)dq\cdot Id}\end{equation}
for $g\in C_0^\infty ({\bf R})$. As in  \cite[Proposition 4.3]{SjZw} 
we can choose
$K(z)$ so that
\begin{equation}
\label{n:7}
{K(z)=e^{i(z-w)Q/h}K(w) \,, }\end{equation}
since this is fulfilled when  \eqref{n:2}, \eqref{n:4} hold (note that
for $ P ( z ) = P - z $, $ U ( z ,w ) = \exp ( i ( z- w ) Q ) $, in 
\cite[(4.4)]{SjZw}).

Consider
\[
{\pi (z)   \stackrel{\rm{def}}{=} {1\over 2\pi h}K(z)K(z)^*.}\] 
Notice that this definition is \indep{} of the choice of $K(z)$ satisfying
\eqref{n:3}.

 We have the properties:
\begin{gather}
\label{n:9}{\pi (z)=\pi (z)^*,}\\
\label{n:10}{P\pi (z)=\pi (z)P=z\pi (z),}\\
\label{n:11}{\int \pi (z)dz=  Id \,, } \end{gather}
where for the understanding of \eqref{n:11}, we remark that if $\phi \in
C_0^\infty ({\bf R})$ is \indep{} of $h$, then $\int\phi (z)\pi (z)dz$ is a
well defined $h$-\pop{} of order 0, and that we can always choose $\phi $
to be equal to 1 in some interval \cor{}ing to the region of phase space
where we wish to establish \eqref{n:11}.

 The relations \eqref{n:9}, \eqref{n:10} are obvious, and so is \eqref{n:11}
 if we work in the
representation \eqref{n:2}, \eqref{n:4}.

We also give a direct proof of \eqref{n:11} which does not involve 
conjugation to the model:
Let $g\in C_0^\infty ({\bf R})$ and consider
\[
{K(z)^*g(Q) K(w),}\] 
where we will assume \eqref{n:7}. Then we get
\begin{equation}
\label{n:13}
{K(z)^*g(Q)K(w)=K^*(z)g(Q)e^{-i(z-w)Q/h}K(z)=({\mathcal F}_hg)(z-w)\cdot Id,}
\end{equation}
where ${\mathcal F}_hg(z)=\int e^{-iqz/h}g(q) dq$ is the $h$-\F{} \tf{}. Here
the second equality follows from \eqref{n:6}.

 Consider the \ml{}ly defined $h$-\pop{}:
$$A=\int \pi (z)dz.$$
$A$ is an elliptic \pop{} of order $0$, and we claim that
\begin{equation}
\label{n:14}
{A^2=A.}\end{equation}
Since the statement is \ml{}, it suffices to show that
\begin{equation}
\label{n:15}
{Ag(Q)A=Ag(Q),}\end{equation}
$g\in C_0^\infty $ (and then choose $g(Q)=1$ in the region of interest).
To get \eqref{n:15}, we first compute:
\[ \begin{split} \pi (z)g(Q)\pi (w)&={1\over (2\pi
h)^2}K(z)K(z)^*g(Q)K(w)K(w)^*\\
&={1\over (2\pi h)^2}K(z)K(z)^*({\mathcal F}_hg)(z-w)e^{i(z-w)Q/h},\end{split} 
\] 
where we used \eqref{n:13}, \eqref{n:6}. Integrating \wrt{} $w$, we get
\[
{\pi (z)g(Q)A=\pi (z)g(Q),} \] 
which gives \eqref{n:14} after integration in $z$.

 Rewriting \eqref{n:14} as $A(A-1)=0$ and using that $A$ is elliptic, we get
\eqref{n:11}.

\medskip
\noindent {\bf Remark.}
 \rm Using \eqref{n:4} (or a direct argument), one can
check that
\[
{\pi (z)={1\over 2\pi h}\int e^{it(P-z)/h}dt,} \] 
where the integral to the right is \ml{}ly well-defined by restricting
the integration to a suitable finite interval.

\medskip
 Let $P,\widetilde{P}$ be two $h$-\pop{}s with simple closed \trajs{}
$$\gamma
\subset p^{-1}(0)\,, \ \ \widetilde{\gamma }\subset \widetilde{p}^{-1}(0)\,\ 
$$ 
with the same period $ T_0 $.
Let $m_0\in \gamma $, $\widetilde{m}_0\in\widetilde{\gamma }$ and assume
that the \cor{}ing monodromy \op{}s 
\[ \begin{split}
& {\mathcal M}(z) \; : \; 
\KKer \longrightarrow \KKer \,, \\
& \widetilde{{\mathcal M}}(z) \; : \; 
\KKert \longrightarrow \KKert \,, \end{split} \]
 are conjugated:

\begin{equation}
\label{n:18}
{{\mathcal U}(z)^{-1}\, \widetilde{{\mathcal M}}(z)\, {\mathcal U}(z)={\mathcal M}(z),}
\end{equation}
where ${\mathcal U}(z): \KKert \longrightarrow \KKer $  is a unitary \fop{}
depending smoothly on $z$.

 Identifying $P,\widetilde{P}$ with $hD_{x_n}$, we see that
\[ {
F  \stackrel{\rm{def}}{=} \int {\mathcal U}(z)\pi (z)dz
}\] 
is an elliptic unitary \fop{} with
\[
{{F_\vert}_{\KKer }={\mathcal U}(z)\,, \ \ \ 
\widetilde{P}F=FP.}\] 

 We can identify $P,\widetilde{P}$ with $hD_{x_n}$ in such a way that
$F$ becomes the identity \op{}. With this identification, we take
$K(z)=\widetilde{K}(z)$ of the form \eqref{n:4}. Without using the identification,
this means that $\widetilde{K}(z)=FK(z)$ and 
\begin{equation}
\label{n:21}
{F={1\over 2 \pi h}\int \widetilde{K}(z)K(z)^* dz.}\end{equation}
Also,
\begin{equation}
\label{n:22}
{\widetilde{M}(z)=M(z).}
\end{equation}

 Let $K_f(z),\widetilde{K}_f(z);K_b(z),\widetilde{K}_b(z)$ be the
forward and backward extensions of $K;\widetilde{K}$ as in \cite[\S 5]{SjZw}.

 If $u$ is \ml{}ly concentrated near $m=\exp (tH_p)(m_0)$, $-\epsilon
\le t\le T_0-2\epsilon $, where $T_0$ is the common period of $\gamma
,\widetilde{\gamma }$, we define $F_fu$, microlocally concentrated to a
\neigh{} of $\exp (tH_{\widetilde{p}}(\widetilde{m}_0))$, by
\begin{equation}
\label{n:23}
{F_fu=e^{-it\widetilde{P}/h}Fe^{itP/h}u.}\end{equation}
Since $F$ intertwines $\widetilde{P}$ and $P$, this definition is invariant
under small variations of $t$ and we see that $F_f$ is a unitary \fop{} with
\canform{}, given by
\[
{\kappa _{F_f}(m)=\exp (tH_{\widetilde{p}})\circ \kappa _F\circ \exp
(-tH_p)(m)\,, }
\]
for $ m $ in a neighbourhood of $ \exp (tH_p)(m_0) $, and 
where $\kappa _F:{\rm
neigh\,}(m_0)\to {\rm neigh\,}(\widetilde{m}_0)$ is the \canform{}
associated to $F$.

 We have
\[
{F_fu={1\over 2\pi h}\int \widetilde{K}_f(z)K_f(z)^*dz.}
\]
In fact, by \eqref{n:21}, \eqref{n:23}:
\[ \begin{split} F_f&={1\over 2\pi h}\int
e^{-it\widetilde{P}/h}\widetilde{K}_f(z)K_f(z)^*e^{itP/h}dz\\   
&={1\over 2\pi h}\int e^{-itz/h}\widetilde{K}_f(z) K_f(z)^*
e^{itz/h}dz={1\over 2\pi h}\int \widetilde{K}_f(z)K_f(z)^* dz.
\end{split}\]

 Similarly, if $u$ is \ml{}ly concentrated near $m=\exp (tH_p)(m_0)$,
$-T_0+2\epsilon \le t\le \epsilon $, we define $F_bu$ to be the \rhs{} of
\eqref{n:23} (now with $t$ in a different interval) and get
\[ {F_bu={1\over 2\pi h}\int \widetilde{K}_b(z) K_b(z)^* dz.}
\] 

 If $m$ belongs to the overlap region, where both the forward and the
backward extensions of $F$ are defined, then for $u$ concentrated to a
\neigh{} of $m$, we have $F_fu=F_bu$:
$$F_fu={1\over 2\pi h}\int \widetilde{K}_f(z)K_f(z)^*udz={1\over 2\pi
h}\int \widetilde{K}_b(z)\widetilde{M}(z) M(z)^*K_b(z)^*udz=F_bu,$$
since
$\widetilde{M}(z)M(z)^*=\widetilde{M}(z)M(z)^{-1}=1$,
in view of \eqref{n:22}.

 This means that $F$ extends to a well-defined \op{} in a neighbourhood of
$\gamma $, and we get

\begin{prop}
\label{p:n.1}
 Let $P,\widetilde{P}$ be two $h$-\pop{}s
with closed simple \trajs{} $\gamma \subset p^{-1}(0)$, $\widetilde{\gamma
}\subset \widetilde{p}^{-1}(0)$ having the same period. 
Let $m_0\in \gamma $, $\widetilde{m}_0\in
\widetilde{\gamma }$, and assume that the monodromy \op{}s
$${\mathcal M}(z) \; : \; \KKer \longrightarrow \KKer $$
$$\widetilde{{\mathcal M}}(z) \; : \; 
\KKert \longrightarrow \KKert $$
are conjugate as in \eqref{n:18}, via an $h$-\fop{} ${\mathcal U}(z)$ depending
smoothly on $z$. Then there exists a unitary \fop{} $F$, \ml{}ly defined
in a \neigh{} of $\gamma $ with an associated \canform{} $\kappa _F:{\rm
neigh\,}(\gamma )\to {\rm neigh\,}(\widetilde{\gamma })$, such that
$ \widetilde{P}F=FP$.
\end{prop} 
\stopthm

Using Theorem \ref{t:1} we will recover the infinitesimal quantum Birkhoff
normal form of $ {\mathcal M } ( z ) $ at $ z = 0$. From the proof of
Proposition \ref{p:n.1} it is clear that the conclusions remain valid to 
infinite order at $ \gamma $ and $ z=0 $, if that assumption is made
on the conjugation. Hence we obtain
\begin{prop}
\label{p:n.2}
Suppose that in Proposition \ref{p:n.1} we only assume that
\[ \widetilde{{\mathcal M}} ( z) \equiv {\mathcal U}(z)\; 
{\mathcal M}(z) \; {\mathcal U}(z) ^{-1} \]  
to infinite order at $m=m_0$, $z=0$ in the sense of \eqref{eq:eq}.
Then there exists a unitary \fop{} $F$, \ml{}ly defined
in a \neigh{} of $\gamma $ with an associated \canform{} $\kappa _F:{\rm
neigh\,}(\gamma )\to {\rm neigh\,}(\widetilde{\gamma })$, such that
\[  \widetilde{P}F=FP \,, \ \ \text{ to infinite order along $ \gamma $}\,.\]
\end{prop}
\stopthm

\section{Applications to inverse problems}
\label{s5}

We start with a simple lemma motivated by \eqref{eq:trace2} and
\eqref{eq:trace3}:

\begin{lem}
\label{l:4.1}
Let $ u ( z, h ) $ be 
a semiclassical family of functions defined near $ z = 0 $,
\[ u ( z, h ) = e^{ i I ( z ) /h } \sum_{j=0}^{\infty} a_j ( z ) h^j 
\,, \ \ I , a_j \in \CI\,, \ a_0 ( 0 ) \neq 0  
\,, \ \Im I(z) = 0\,. \]
Suppose that $ {\mathcal V} \subset \CIc ( \RR ) $ is a subspace
from which we can pick an element having any prescribed finite Taylor
polynomial at $ I'( 0 ) $. If for any $ g \in {\mathcal V} $ we know
$ b_j (  g ) $ where
\[  J ( g, u ) = h^{-1} \int \hat g ( z / h ) \chi ( z ) u ( z , h ) dz =
e^{ i I ( 0 ) / h } \sum_{k = 0 }^\infty b_k ( g ) h^k \,, \ \
\chi \in \CIc ( \RR ) \,, \ \chi \equiv 1 \ \text{near $ 0$,} \] 
then the distributions $ b_k $ determine the $ I^{ ( l )} ( 0 ) $,
and $ a_j ^{(l)} ( 0 ) $, for all $l $ and $ j $.
\end{lem}
\begin{proof}
We first observe that if $ a_j ( z ) = {\mathcal O} ( z^\infty ) $ then
they contribute $ {\mathcal O} ( h^\infty ) $. Hence in considering 
the expansion of $ J ( g , u ) $ we can replace $ I $ and $ a_j $'s by 
their formal power expansions:
\[  I ( z ) = I_0 + I_1 z + \sum_{ k = 1}^\infty  I_{k+1} z^{k+1} \,, \ \
a_j ( z)  = \sum_{k = 0}^\infty a_{ j k } z^k \,, \]
and our task is the recovery of $ I_k $'s and $ a_{ j k } $'s.
Making a change of variables $ \zeta = z / h $ we write $ J ( g , u ) $
as the formal integral
\[ e^{ I_0 / h } \int_{\RR} \hat g  ( \zeta ) e^{ i I_1 \zeta }
e^{ i \sum_{  j= 1}^\infty h^j I_{ j + 1} \zeta^{ j + 1} }
\left( \sum_{ j, k \geq 0 } a_{ jk} \zeta^k h^{ k + j } \right)
\; d \zeta \,,\]
where we can neglect $ \chi ( h \zeta ) $, as $ \chi ( z ) \equiv 
1 $ near $ z = 0 $.
Just as in the proof of Theorem \ref{t:1} we now see that the coefficient of
$ h^{p+1} $ in the expansion is given by
\[ e^{ I_0/ h } h  \int \hat g ( \zeta ) e^{ i I_1 \zeta } 
\left( \sum_{ l + k = p }  a_{ lk } \zeta^l + 
i a_{00} I_{p + 1} \zeta^{p+1} + H_m ( a , I ) \right) \,, \]
where $ H_m  ( a , I ) $ depends only on $ a_{ lk} $ with $ l + k < p $,
and $ I_{ l} $ with $ l < p+1 $. Since we assumed that $ a_{00} = a_0 ( 0 ) 
\neq 0 $, and that we can find $ g \in {\mathcal V} $ with any 
prescribed finite Taylor polynomial
at $ I_1 $, the coeffients of $ I $ and $ a $
can be recovered from $ b_j  ( g) $'s, $ g \in {\mathcal V} $.
\end{proof} 

Using Lemma \ref{l:4.1}, and \cite[Theorem 2, Proposition 7.5]{SjZw} 
(see \eqref{eq:trace2} \eqref{eq:trace3} above) we obtain the
following consequence of Theorem \ref{t:1}: 

\begin{thm}
\label{t:4}
Suppose that $ P ( z ) $ satisfies the assumptions needed for 
\eqref{eq:trace2} and is classical in the sense of \eqref{eq:clas}. 
Suppose that in addition the eigenvalues of the linear Poincar{\'e} map 
of $ \gamma ( 0 ) $ satisfy \eqref{eq:dio2}. Then the coefficients in the
asymptotic expansion of 
\[ \frac1\pi  {\mathrm{tr}} \;  \int f ( z / h ) 
 \bar \partial_z \left[ \tilde \chi ( z )\;   
\partial_z P (z) \;  P ( z) ^{-1} \right] A 
{\mathcal L} ( d z )  = 
\sum_{ k \in {\mathbb Z} \setminus \{ 0 \} } 
 e^{ i k I  / h } \sum_{j = 0}^\infty a_{j, k} ( f) h^j  \,,
\]
determine the {\em infinitesimal quantum Birkhoff normal form}  of the 
quantum monodromy operator $ M ( z ) $ at $ z = 0 $. 
In particular the Birkhoff normal form of the Poincar{\'e} map
 of $ \gamma ( 0 )$ is determined.
Here the notation is the same as in \eqref{eq:trace2}, and
by the coefficients we mean the distributions, $ f \mapsto a_j ( f ) $, 
$ \hat f \in \CIc ( \RR \setminus \{ 0 \} ) $.
\end{thm}
\begin{proof}
We first see that the non-degeneracy of $ \gamma ( z ) $ for $ z $ close
to $ 0 $ shows that \eqref{eq:trace3} holds.
We apply \eqref{eq:trace2} and integrate by parts:
\[ \begin{split}  
& - \frac{1 }{ 2\pi i } 
\sum_{
-1\ne
 k=  -N -1  }^{N -1 }  {\mathrm {tr}} \; 
\int_{\mathbb R} f (z/h) M(z, h)^k \frac{d}{dz} M ( z , h )
 \chi ( z) 
dz = \\
& \ \ \ \ \ \ h^{-1}
\frac{1 }{ 2\pi i } 
\sum_{
 -1\ne 
k = -N -1  }^{N -1 }  {\mathrm {tr}} \; 
\int_{\mathbb R} f' (z/h) (k+1)^{-1} M(z, h)^{k+1} 
 \chi ( z)  + {\mathcal O} ( h^\infty )  \,, \end{split} \]
where $ \chi \equiv 1 $ near $ 0$ and hence the contribution of $ \chi' $ is
negligible: $ f ( z / h ) = {\mathcal O} ( h^\infty ) $ for 
$ |z| > \epsilon > 0 $. The assumptions of Lemma \ref{l:4.1} with 
$ \hat g = f' $ are satisfied and hence  we obtain expansions of
$ \tr M ( z , h )^k $, for any $ k $.
Now we can apply Theorem \ref{t:1} with $ U(z) = M ( z , h ) $ and the 
proof is completed.
\end{proof}

When we combine this result with Proposition \ref{p:n.2} we obtain 
infinitesimal information about $ P ( z) $ at $ \gamma$:

\begin{coro}
\label{c:n}
Suppose that the assumptions of Theorem \ref{t:4}
are satisfied and $ P ( z) = A ( z)^* ( P - z ) A ( z)$, where $ A ( z ) $ is 
elliptic near $ \gamma $. Then, up to conjugation by elliptic $h$-Fourier
integral operators, in the sense of Proposition \ref{p:n.2}, $ P $ 
is determined to infinite order at $ \gamma $.
\end{coro}

When specialized to the classical, or high energy, case, we obtain
\begin{coro}
\label{c:new}
Let $ M $ be a compact $ \CI $ manifold 
and let $ P $ be a positive, self-adjoint
elliptic classical pseudo-differential operator acting on  $ \CI ( 
M , \Omega_M^{\frac12}) $. Suppose that
$ \gamma  \subset p^{-1} ( 1 ) $ is a simple non-degenerate orbit of the classical flow of the principal 
symbol of $ P $, and $ P_\gamma $ is its linear Poincar{\'e} map.  If
\eqref{eq:dio2} holds for the
eigenvalues of $P_\gamma$ \eqref{eq:ev}, and no other orbit in 
$ p ^{-1} (1) $ has the period
which is a multiple of the period of $ \gamma $, 
then the spectrum of $ P $ determines the Birkhoff normal form of the
Poincar{\'e} {\em first return} map. 
\end{coro}
\begin{proof} 
Suppose that the order of $ P $ is $ m$. We define a
{\em semi-classical operator}, $ P ( h ) $, by putting,
$ P ( h ) = h^{  m } P ( x ,  D_x ) -1 \,, $
so that the (semi-classical) 
principal symbol of $ P ( h ) $ is the same as the (classical) principal
symbol of $ P $. The non-degeneracy assumption for $ \gamma \subset 
\{ p = 1\} $ implies that the assumptions of Theorem \ref{t:4} are 
satisfied, with $ P ( z) = P ( h ) - z $. The simplicity of $ \gamma $ then
shows that for $ f $'s with $ {\rm{supp} }\; \hat f $ compact, and 
sufficiently close to $ \ZZ T_\gamma  $, where $ T_\gamma $ is the
(primitive) period of $ \gamma $ on $ \{ p = 1 \} $, we have
\begin{equation}
\label{eq:loc} \tr \; f ( P ( h ) / h ) \chi ( P ( h ) ) A = 
 \tr \; f ( P ( h ) / h ) \chi ( P ( h ) )  + {\mathcal O} ( h^\infty ) \,.
\end{equation}
The left hand side is a (semi-classical) spectral invariant and hence
the result follows from Theorem \ref{t:4}.
\end{proof}

We leave it to the reader to formulate the corresponding corollary in the more
general semi-classical case.

In the case of manifolds with boundaries we consider, for the sake of
simplicity, the case of second order differential operators only:
\begin{thm}
\label{t:3}
Suppose that $ (M, g ) $ 
is a Riemannian manifold with boundary, and that $ \Delta_g $ is
its Dirichlet or Neumann Laplacian.

Let 
$ \gamma $ be a simple non-degenerate orbit of the classical billiard 
flow of $ g$.
If all the intersections of $ \gamma $ with $ \partial M $ are {\em 
transversal},  \eqref{eq:dio2} holds for 
eigenvalues of the Poincar{\'e} map of $ \gamma$ \eqref{eq:ev}, 
and the lengths of other closed orbits of the generalized billiard 
flow do not accumulate
at a multiple of the period of $ \gamma $, 
then the spectrum of $  \Delta_g  $ determines the Birkhoff normal form of the
Poincar{\'e} {\em first return} map of $ \gamma $.  
\end{thm}

\noindent
{\em Outline of the proof:}
First we need to indicate why the results of \cite{SjZw} 
can be applied in this situation.
For that we will check that local solutions can be
microlocally defined near a trajectory $ \gamma $.  Let $ p $ 
be the symbol of $ D_t^2 - \Delta_g $ and $ q $ the defining function of 
$ \partial M $. The  transversality condition
 on $ \gamma $ means that $ H_q p \neq 0 $ at the boundary. We first
recall (see for instance \cite[Sect.24.2]{Horb}) that we can find coordinates
in which, near a point in $ \partial M $, $ q = x_1 $, and 
$ p ( x, \xi ) = \xi_1^2 - r ( x , \xi' )$, $ \xi' = ( \xi_2, \cdots, 
\xi_{m +1 } )$, $ m = {\rm dim}\; M $. 
Solving the boundary problem microlocally is equivalent to a hyperbolic problem and hence
we have a solution. For localization needed in \eqref{eq:trace2} 
we can use, near the boundary, an operator $ A = A ( x, h D_{x' } ) $.
\footnote{Alternatively, we could have reduced the problem to the boundary,
and considered the monodromy operator there.}

A more subtle argument is needed for showing that \eqref{eq:loc} holds,
as that involves results on propagation of singularities for boundary
value problems -- see \cite[Sect.24.5]{Horb} and references given there.
These results provide the foundation for 
 the trace formula of Guillemin and Melrose \cite{GM} (see 
also \cite{PS}), which 
among other things shows that $ {\rm{singsupp}}\,
 ( \tr  \cos ( t \sqrt {-\Delta_g} )  ) \subset 
\overline{\mathcal L } $, where $ {\mathcal L } $ is the set of length of
closed orbits of the billiard flow. It also follows from \cite{GM}
and \cite{PS} that if 
$ A ( x , \langle D_t \rangle^{-1} D_x ) $ is an operator which provides
a localization to a neighbourhood of the bicharacteristic strip of $ \gamma $,
then 
\begin{equation}
\label{eq:gm}
 \tr \cos ( t \sqrt { - \Delta_g} ) - \tr \cos ( t \sqrt {- \Delta_g} ) 
 A ( x , \langle D_t \rangle^{-1} D_x )  \in \CI ( \RR ) \,, 
 \ \ \text{ for $ t $ near 
$ k T_\gamma $,} \end{equation}
that is, 
the singularities of the trace of $ \cos ( t \sqrt {-\Delta_g} ) $ 
at $ k T_\gamma $
come from a small neighbourhood of $ \gamma $ (depending on $ k $).
As stated above we only need tangential pseudo-differential operators
to localize near $ \gamma$.

A straightforward translation to the semi-classical
setting (see for instance \cite{Zw})  gives \eqref{eq:loc}. In fact, by 
standard oscillatory testing,
$  A ( x, \langle D_t \rangle^{-1} D_x ) e^{ i t / h } = 
 e^{ it / h } \widetilde A ( x , h D_x ; h )  \,, $
where $ \widetilde A $ has the  semi-classical principal symbol given by 
$ A ( x , \xi ) $.
This shows that for $ \hat f \in \CIc ( \RR ) $ even, 
\[ \tr f ( ( h \sqrt {- \Delta_g}  - 1 ) /h ) 
 \tilde A ( x , h D_x , h) =
\frac{1}{2 \pi } \int \tr \left(  A ( 
 x, \langle D_t \rangle^{-1} D_x ) \cos t \sqrt {- \Delta_g} \right) 
\hat f ( t ) e^{ i t / h } dt \,,\]
and thus \eqref{eq:gm} follows from \eqref{eq:loc}. 
\stopthm

\medskip

By combining Theorem \ref{t:4}, or rather Corollary \ref{c:new},
with Proposition \ref{p:n.2} or with 
the geometric arguments of \cite[\S 2]{G} we obtain
the {\em Weinstein conjecture} first proved by Guillemin
\cite{G} in the case of elliptic trajectories, and by Zelditch 
\cite{Zel1},\cite{Zel2} in general. It involves the concept of
the {\em Birkhoff normal form for the Hamiltonian} -- 
\cite[Theorem 1.1]{G}. As shown 
by Fran{\c c}oise and Guillemin \cite{GF} and Guillemin \cite[\S 2]{G},
and as can be also deduced from Sect.\ref{s4} above,
that normal form is essentially equivalent to the normal form 
of the Poincar{\'e} map.

\begin{coro}
\label{c:1}
Let $( M, g ) $ be a compact Riemannian manifold without boundary, 
and $ \gamma $ a simple
non-degenerate geodesic on $ M $, with the spectrum of its linear
Poincar{\'e} map satisfying \eqref{eq:dio2}.
Then the spectrum of the Laplacian, 
$ \Delta_g $, determines the Birkhoff normal form of $ \gamma $.
\end{coro}

In the special case of symmetric domains in the plane, we can apply
Theorem \ref{t:3} and  
an observation of Colin de Verdi{\`e}re \cite[\S 4]{CdV} to recover a result of
Zelditch \cite{Zel}\footnote{Melrose \cite{Mel} has announced the same
result, sketching a yet different proof. Recently, by considering the full 
quantum Birkhoff normal form, Zelditch \cite{Zel-p}
succeeded in removing one of the
symmetry conditions}:

\begin{coro}
\label{c:2}
Suppose that $ \Omega \subset \RR^2 $ is a bounded domain,
$ \partial \Omega $ is real analytic, and that $ \Omega $ is 
symmetric with respect to the $ x $ and $ y $ axes.
Suppose also that $ \Omega $ has a simple non-degenerate bouncing ball orbit
along one of the axes. 
Then the spectrum of
the Dirichlet (or Neumann) Laplacian determines $ \Omega $ among all
such domains.
\end{coro}

\end{document}